\input amstex
\documentstyle{amsppt}
\NoRunningHeads 

\magnification \magstep 1 
\def\r{\roman{r}}
\def\Zt{\Bbb Z}
\def\Pt{\Bbb P}
\def\Rt{\Bbb R}
\def\Ct{\Bbb C}
\def\Nt{\Bbb N}
\def\Qt{\Bbb Q}
\def\Rp{\accent'137Z}
\def\lp{\char'40l}
\def\ap{a\kern -.32em\lower -0.1ex\hbox{\char'030}\kern -.4em\phantom{.}}

\topmatter
\title Residue in intersection homology for quasihomogeneous 
singularities\endtitle
\author Andrzej Weber \endauthor
\thanks Partially supported by KBN 2P30101007\endthanks
\affil Institute of Mathematics, University of Warsaw \endaffil
\address ul. Banacha 2, 02--097 Warszawa, Poland \endaddress
\email aweber\@mimuw.edu.pl \endemail
\date  October 1996\enddate
\endtopmatter
\document

\head 1. Residues in homology and the problem of lift \endhead

Let $M^{n+1}$ be a complex manifold and let $K$ be a hypersurface in $M$. 
Suppose
$K$ has isolated singularities. We think of $K$ as obtained from a
manifold--with--boundary $K^\circ$ by shrinking the components of the 
boundary.
For simplicity assume that $n>1$. Then the intersection homology of $K$ is
isomorphic to \cite{7}, \cite{4}:
$$\aligned IH^{\underline m}_n(K)&= im\left([K]\cap:H^n(K)@>>>H_n(K)\right)=\\
&=im \left([K^\circ]\cap:H^n(K^\circ,\partial K^\circ)
@>>>H_n(K^\circ,\partial K^\circ)\right) =\\
&=im \left(H^n(K^\circ,\partial K^\circ) @>>>H^n(K^\circ)\right) =\\
&= ker\left(H^n(K^\circ)\longrightarrow H^n(\partial K^\circ)\right).
\endaligned $$
In \cite{13} we define a residue of a closed form
$\omega\in\Omega^{n+1}(M\setminus K)$ as the class of the Leray \cite{9} residue 
form
$res\,\omega=[Res\,\omega]\in H^n(K\setminus\Sigma)\simeq H_n(K)$.
All coefficients of homology and cohomology are in $\Ct$. 
The following question arises: can one integrate the residue form over cycles
intersecting singularity. The question comes from partial differential
equations; \cite{14}. It turns out, that certain integrals have a meaning, but it is
not clear in which spaces cycles and integrals should be considered.
The possibility to give a meaning to the symbol $\int_\xi Res\, \omega$ which
would not depend on the homology class of $\xi$ is simply a lift of $[Res\,
\omega]$ to cohomology:
$$\def\awr{\mathrel{\smash -}}
\def\przearrow{\awr\awr\rightarrow}
\matrix H^n(K\setminus\Sigma) & \simeq &H_n(K) & @<\;PD\;<<& H^n(K)\\
[Res\, \omega]&=&res\,\omega &\przearrow& ?\endmatrix$$
The homological residue $res\,\omega$ can be defined equivalently by the
Alexander duality: 
$$\matrix H^{p+1}(M\setminus K) &@>\;\;\delta\;\;>>&H^{p+2}(M,M\setminus K)& @>
[M]\cap>>&H_{2n-p}(K)\\
\omega &\longmapsto & \delta\omega & \longmapsto & res\,\omega \endmatrix\,.$$

The previous case was $p=n$. Now $p$ is arbitrary and $K$ is a variety with
not necessary isolated singularities.
The question is the same: does the homological residue $res\, \omega $ belong 
to
the image of Poincar\'e duality map. The Poincar\'e duality map factors
through the intersection homology:
$$
\def\mpr#1{\smash{\mathop{\hbox to 25pt{\rightarrowfill}}\limits^{#1}}} 
\matrix H^{2n-p}(K) &\mpr{PD} & H_p(K) \\
 \searrow&&\nearrow\\
&IH^{\underline m}_p(X)\,.& \endmatrix$$
The question of lift to intersection homology also seems to be 
reasonable. In many
cases algebraically defined elements in $H_*(K)$ lift to $IH^{\underline 
m}_*(K)$, e.g.
Chern-MacPherson classes or arbitrary algebraic cycles; \cite{3}. 

Consider a case
when $M$ is a projective manifold and
$K=\bigcup_{i\in I}D_i$ is a sum of smooth divisors with normal 
crossings. Then
$$IH^{\underline m}_p(K)=\bigoplus_{i\in I}H^{2n-p}(D_i)\,.$$
One can easily show that:

\proclaim{Proposition 1 }Homological residue of a class $c\in H^p(M\setminus K)$
can be lift to the intersection homology if and only if $c$ belongs to 
 $W_{p+1}H^p(M\setminus K)$ -- the
$(p+1)$--st term of the Deligne weight filtration. \endproclaim

Belonging to $W_{p+1}H^p(M\setminus K)$ means that $c$ is represented by a 
smooth form
$$\omega \in \Omega^{p-1}(M)\wedge \Omega^1(log <K>)\,,$$ 
see e.g. \cite{8 \S 5}. The above proposition
remains true for an arbitrary singular projective va\-rie\-ty. To prove this 
one can 
use a
resolution of $K$, then push the residue to the intersection homology of the
resolution and pull down using the result of \cite{3}. Unfortunately this
procedure uses desingularization, weight filtration and functoriality of
intersection homology. Each of these ingredients is rather 
mysterious and hard
to compute. We will not follow this direction. We restrict our attention to
the case of isolated quasihomogeneous singularities.

I would like to thank
several people who have influenced me during the work on residues. First of
all I would like to thank Professors B. Ziemian and H. \Rp o\lp \ap dek,
P.Jaworski, J-P. Brasselet, G. Barthel.

\head 2. Valuation and quasihomogeneous functions \endhead

Let $\Ct(z_0,\dots ,z_n)$ be the field of the rational functions on $n+1$
variables and let $$v:\Ct(z_0,\dots ,z_n)@>>>\Qt$$ be a valuation satisfying
\item{1)} $v(z_i)=a_i>0$;
\item{2)} $v(f\,g)=v(f)\,v(g)$;
\item{3)} if $f=\sum f_i$, $f_i$ monomial, then $v(f)=max\{\,v(f_i)\,\}$

If $f$ is a sum of monomials of the same weight, then we say that $f$ is
quasihomogeneous (with respect to the valuation $v$).
We define the valuation on forms 
$$v:\Ct(z_0,\dots ,z_n)\otimes\Lambda^*\left((\Ct^{n+1})^*\right)@>>>\Qt$$ 
putting $v(dz_i)=v(z_i)=a_i$.

Suppose that in each singular point the hypersurface $K$ is given by an 
equation $s=0$
with $s$ quasihomogeneous (in some coordinates and valuation). We assume 
that
$v(s)=1$. We define a number 
$$\kappa=v(dz_0\wedge \dots \wedge dz_n) = \sum_{i=0}^n a_i\,.$$
It coincides (after the change of the sign) with 'the complex oscillation 
indicator'
defined in \cite{1; 13.1 p.258} and it is an analytic invariant of a 
singularity type.
We define a condition:

\proclaim{Condition 2} For any choice of $k_i\in\Nt\cup \{0\}$, $i=0,\dots ,n$
we have $\kappa+\sum k_ia_i\not=1$.\endproclaim

Of course the Condition 2 is satisfied if $\kappa>1$.

\example{Example 3} Let $s=z_0^3+z_1^3+z_2^4$. Then $a_0=a_1=\frac 13$ 
and
$a_2=\frac 14$, $\kappa=\frac {11}{12}$. The Condition 2 is satisfied.\endexample

\head 3. A simple criterion to lift \endhead

We investigate residues of meromorphic forms of the type $(n+1,0)$ with a
first order pole 
on $K$, i.e. the forms which can be locally written as $\omega=\frac gs
dz_0\wedge\dots \wedge dz_n$ with $g$ holomorphic.

\proclaim{Theorem 4} Suppose that $K$ of dimension $n$ has isolated
singularities given by 
quasihomogeneous equations. Let $\omega\in\Omega^{n+1,0}(M\setminus K)$ be a
meromorphic form with a first order pole on $K$.
If the Condition 2 is fulfilled in each singular
point then the residue class of $\omega$
lifts to intersection homology of $K$.\endproclaim

The Theorem 4 is true in a greater generality; we assume that $K$ has isolated
singularities and the Condition 2 should be substituted by the following:

\proclaim{Condition 5} The number 0 does not belong to the spectrum of
the singularity. \endproclaim

The concept of the spectrum as defined in \cite{1; \S 13.3 p. 270} comes from the 
theory of
oscillating integrals. It follows from the definition that
$\int_XRes\,\omega=0$ for any cycle contained in the link of a singular point;
see \cite{13; 1.3}. In this case residue lifts to the cohomology. The point is
that the spectrum is computable from the Newton diagram; \cite{1; \S 13.3 p. 274}. 
For
quasihomogeneous singularities we have 
$$\gather \{ \text{spectrum of the singularity}\} \cap (-\infty ,0] =\\
=\{ v(\frac gs dz_0\wedge\dots\wedge dz_n)\} \cap (-\infty ,0] =\\
=\{ \kappa +\sum k_ia_i-1 : \, k_i\in\Nt\cup \{0\},\; i=0,\dots ,n\} \cap
(-\infty ,0] \endgather$$ 
The spectrum can also be defined in terms of eigenvalues of the monodromy
acting on the vanishing cycles filtered by weights; \cite{11}.

Now we prove the Theorem 4 using few well known facts from the intersection
homo\-logy theory.

\demo{Proof} One should show that $[Res\,\omega]\in
ker\left(H^n(K^\circ)\longrightarrow H^n(\partial K^\circ)\right)$ that is 
for each link $L$ in $K$
$[Res\,\omega_{|L}]=0\in H^n(L)$. We take a neighbourhood
of a singular point in which $K$ is given by a quasihomogeneous equation
$s=0$. We give a formula for $Res\,\omega$ in the points where
$s_0=\frac{\partial s}{\partial z_0}\not=0$: 
$$\aligned ds&=\sum_{i=0}^n s_idz_i\,,\\
dz_0&=\frac 1{s_0}ds-\sum_{i=1}^n\frac{s_i}{s_0}dz_i \\
\omega &=\frac gs \frac{ds}{s_0}\wedge dz_1\wedge\dots\wedge dz_n=\\
&={ds\over s}\wedge {g\over s_0} dz_1\wedge\dots\wedge dz_n\,.\endaligned$$

Let $\r=\frac g{s_0} dz_1\wedge\dots\wedge dz_n$. Then $Res\,\omega=\r_{|K}$
in the points where $s_0\not=0$. Now suppose that $g$ is quasihomogeneous (or
we decompose $g$ into a quasihomogeneous components). Let $v(g)=\alpha$.
We have
$$v(\omega)=v(ds)-v(s)+v(\r)=v(\r)\,.$$
Then $$\aligned 
v(\r )=&v(g)-v(s)+v(dz_0)+\dots +v(dz_n)=\\
=&\alpha-1+a_0+\dots +a_n= \alpha-1+\kappa\,. 
\endaligned $$

Let $l $ be a natural number such that $l \,a_i\in \Nt$ for $i=0,\dots , n$. 
We
construct a branched covering of $\Ct^{n+1}$:

$$\aligned \Phi: \Ct^{n+1}&\longrightarrow\Ct^{n+1}\\
\hat z_0, \dots ,\hat z_n&\longmapsto \hat z_0^{l a_0}, \dots ,\hat
z_n^{l a_n}\,.\endaligned $$
Let $\hat v$ be a standard valuation: $\hat v(\hat z_i)=1$. The map
$\Phi$ has the property: $$\hat v(\Phi^*\eta)=l \,v(\eta)$$ for any $\eta\in
\Ct(\hat z_0,\dots ,\hat z_n)\otimes\Lambda^*\left((\Ct^{n+1})^*\right)$.
We have $$v(\Phi^*\r)=l (\alpha-1+\kappa)\,.$$
If we write $\Phi^*\r=q\,d\hat z_1\wedge\dots\wedge d\hat z_n$ then $q$ is
homogeneous function of weight $$v(q)=l (\alpha-1+\kappa)-n\,.$$
The mapping $\Phi$ induces the branched covering of the links:
$$ S^{2n+1}\cap\Phi^{-1}(K)=\widehat{L} @>\bar \Phi>>
L=K\cap \{z_0,\dots , z_n: |z_0|^{2a_0}+\dots +|z_n|^{2a_n}=1\}$$
The degree of this map is $l \,\kappa$. Unfortunately $\widehat{L}$ may
be singular; see Example 6. To show that $[Res\,\omega_{|L}]=0$
we will prove that $ [\bar \Phi^* Res\,\omega_{|L}]=0\in
IH^{\underline m}_{n-1}(\widehat{L})$.
It is enough since the map
$$H^n(L)@>\bar \Phi^*>>H^n(\widehat{L})@>>>IH^{\underline
m}_{n-1}(\widehat{L})$$ 
is a monomorphism with a splitting
$$IH^{\underline m}_{n-1}(\widehat{L})@>>>H_{n-1}(\widehat{L})@>\bar \Phi_*>>
H_{n-1}(L)@>>>H^n(L)\,.$$ 
The last map is the inverse to the Poincar\'e duality isomorphism
multiplied by $(l \,\kappa)^{-1}$, the maps to and from
intersection homology are the canonical ones. To show vanishing
in intersection homology we use a Gysin sequence of the
fibration 
$$S^1\hookrightarrow\widehat{L}@>p>>\widehat{L}/S^1$$
coming from the action of $\Ct^*$ on $\Phi^{-1}(K)$:
$$@>>>IH^{\underline m}_{n}(\widehat{L}/S^1)@>\cap e>>
IH^{\underline m}_{n-2}(\widehat{L}/S^1)@>p^*>>
IH^{\underline m}_{n-1}(\widehat{L})@>p_*>>
IH^{\underline m}_{n-1}(\widehat{L}/S^1)@>>>\,.$$
The map $\cap e$ is the multiplication by the Euler class of the
fibration; it is an isomorphism by hard Lefschetz since
$\dim_\Ct\widehat{L}/S^1 = n-1 $; \cite{2}. We interpret $IH^{\underline
m}_{n-1}(\widehat{L})$ 
as the $L_2$--cohomology of the nonsingular part of $\widehat{L}$:
$$IH^{\underline m}_{n-1}(\widehat{L})=H^n_{(2)}(\widehat{L}\setminus\Sigma)
=:H^n_{(2)}(\widehat{L})$$
for suitably chosen metrics on $\widehat{L}\setminus\Sigma$ and $(\widehat 
{L}\setminus\Sigma)/S^1$; see \cite{6}, \cite{12}. Then the sequence has a form:
$$@>>>H_{(2)}^{n-2}(\widehat{L}/S^1)@>\cap e>>
H_{(2)}^n(\widehat{L}/S^1)@>p^*>>
H_{(2)}^n(\widehat{L})@>p_*>>
H_{(2)}^{n-1}(\widehat{L}/S^1)@>>>\,.$$
The map $p_*$ is just the
integration along the fibers of $p$. Let us calculate the
integral in the trivialization of the bundle $\Ct^{n+1}\setminus \{0\} @>p>>
\Pt^n $ over $ U_0=\{\hat z_0\not=0\} \subset \Pt^n$:
$$\aligned
\Ct^*\times U_0 &\longrightarrow p^{-1}(U_0)\,,\\
u_0, u_1,\dots ,u_n &\longmapsto u_0, u_0u_1,\dots
,u_0u_n\,.\endaligned$$ 
We write $\Phi^*\r$ in $u$--coordinates:
$$\aligned \Phi^*\r&=q(\hat z_0,\dots ,\hat z_n)d\hat 
z_1\wedge\dots\wedge d\hat 
z_n=\\
&= u_0^{l (\alpha-1+\kappa)-n}\bar q(u_1,\dots ,u_n)
(u_1du_0+u_0du_1)\wedge\dots\wedge (u_ndu_0+u_0du_n)= \\
&= u_0^{l (\alpha-1+\kappa)-1}\bar q(u_1,\dots ,u_n)
du_0\sum_{i=1}^n(-1)^{i+1}u_idu_1\wedge\dotsb\overset i\to\vee \dotsb\wedge
du_n +\\
&\hskip 239pt+u_0du_1\wedge\dots\wedge du_n=\\
&=u_0^{l (\alpha-1+\kappa)-1}du_0\wedge\r_2+\Theta\,,\endaligned$$
where $\r_2$ and $\Theta$ do not contain $du_0$ and $\r_2$ does not
depend on $u_0$. Thus the integral can be nonzero only if
$\alpha+\kappa=1$. It is impossible by the Condition 2 since
$\alpha$ is a combination of $a_i$'s. Thus
$p_*\Phi^*(Res\,\omega_{|L})=0$, so the residue lifts to intersection
homology. \enddemo

\example{Example 6} Consider the polynomial
$$s(x,y,z)=(x+z^2)^2+y^2-z^4\,.$$
It has an isolated singularity of the type $A_3$. It is
quasihomogeneous with weights $v(x)=v(y)=\frac 12$ and
$v(z)=\frac 14$. The polynomial $\Phi^*(s)$ is:
$$\Phi^*(s)=(x^2+z^2)^2+y^4-z^4=x^4+2x^2z^2+y^4\,.$$
Zero is not an isolated singularity since for $z=c=\text{const}$ we obtain:
$$x^4+2x^2c^2+y^4 \sim x^2+y^4$$
which is a singularity of the type $A_3$.
The example of a singularity with $\widehat{L}$ nonsingular is 
$z_0^{k_0}+\dots +z_n^{k_n}$ for any choice of $k_i\in \Nt$.\endexample

The Example 6 shows, that in the proof of the Theorem 3 we have
to use the hard Lefschetz theorem for intersection homology
instead of the standard one.

\remark{Remark} Note that (in the case of isolated
singularities) if the residue class lifts to intersection homology
then it lifts to cohomology. Choose $p>1$. For $\kappa>1$ we
show that $Res\,\omega$ is a form with the $L_p$--integrable
norm in suitably chosen metric; \cite{13}. If $p$ is large then the
$L_p$--cohomology of $K\setminus\Sigma$ is isomorphic to the
cohomology $K$; see \cite{5} ,\cite{12}. This way we find a particular lift
to cohomology. This lift depends on coordinates, but can it be
calculated in terms of integrals. \endremark

\head 4. Weighted blow-up \endhead

There is another way of looking at the calculation presented in the proof of
the Theorem 4. Let the group $G=\Zt/la_0\times\dots\times \Zt/la_n$ acts on 
the
coordinates of $\Ct^{n+1}$ by the multiplication by the roots of unity. Then
$K=\widehat{K}/G$. We blow up $\widehat{K} \subset 
\Ct^{n+1}$ in $0$ and obtain:
$$\matrix \widetilde{\Ct^{n+1}}&\supset &\widehat{Y}\cup \Pt^n&
@>\widetilde{\Phi}>> &Y\cup \Pt(v) & =& 
\widehat{Y}/G\cup \Pt^n/G&\subset &\widetilde{\Ct^{n+1}}/G\\ 
@V\widehat{pr}VV @VVV @VVV && @VprVV\\
\Ct^{n+1}&\supset &\widehat{K} &@>\Phi>> &K&=& \widehat{K} /G&\subset & 
\Ct^{n+1}/G&
=&\Ct^{n+1}. \endmatrix$$
Here $\Pt(v)=\Pt^n/G$ is weighted projective space. We have $\widehat{Y}\cap
\Pt^n= \widehat{L}/S^1$ and $Y\cap\Pt (v)= L/S^1$.
The spaces $\Pt(v)$, $\widetilde{\Ct^{n+1}}/G$, $Y$ and $L/S^1$ are
homology manifolds; locally they are quotients of smooth manifolds by a finite
group i.e. they are V--manifolds as defined by Steenbrink; \cite{10}. From the
homology point of view they can be treated as ordinary (smooth) K\"ahler
manifolds. 

\head 5. Nonvanishing of the second residue \endhead

The last lines of the proof of the Theorem 4 lead to a definition of an
element $$res_2\omega=\left[\frac 1{2\pi i}\int_p Res\,\omega_{|L}\right]\in
IH^{\underline m}_{n-1}(\widehat{L}/S^1)\,.$$
This is an obstruction to lift the residue class to cohomology. We call it the
second residue. The class $res_2\omega$ is $G$--invariant, so it is in 
$$IH^{\underline m}_{n-1}(\widehat{L}/S^1)^G=IH^{\underline 
m}_{n-1}(L/S^1)=H^{n-1}(L/S^1)\,.$$ The form
$r_2$ represents 
$res_2\omega$.
Since $L/S^1$ is V--manifold then its cohomo\-logy admits Hodge decomposition 
\cite{10}
and $res_2\omega$ is of $(n-1,0)$ type. The form $r_2$ is harmonic outside 
the 
singularities of $L/S^1$, so to show that it does not vanish in cohomology it
suffices to check that it is not tautologically zero.
Suppose that $g$ is quasihomogeneous of the weight
$$v(g)=\alpha=1-\kappa\,,$$
hence $v(\omega)=0$. We will show that 
$$res_2\omega=res_2\left(\frac gsdz_0\wedge\dots\wedge dz_n\right)\neq 0 
\in H^{n-1}(L/S^1)\,.$$ 

We blow--up $\Ct^{n+1}\supset \widehat{K}$ in $0$. We calculate the form 
$\Phi^*\omega$
pulled up to $\widetilde{\Ct^{n+1}}$ in the canonical coordinates (in the 
$0$--th
chart). 
$$\gather 
\widehat{pr}^*\Phi^*\omega=C\frac{\Phi^*g}{\Phi^*s}\left(\prod_{i=0}^{i=n}\hat
z_i^{la_i-1}\right)\, d\hat z_0\wedge\dots\wedge d\hat z_n=\\
=C\frac{u_0^{l\alpha}\widetilde{\Phi^*g}}
{u_0^l\widetilde{\Phi^*s}}u_0^{l\kappa-n-1}\left(\prod_{i=1}^{i=n}
u_i^{la_i-1}\right)\, u_0^n\,du_0\wedge\dots\wedge du_n=\\
=C\frac{du_0}{u_0}\wedge\frac{\widetilde{\Phi^*g}}{\widetilde{\Phi^*s}}\left(
\prod_{i=1}^{i=n}u_i^{la_i-1}\right)\, du_1\wedge\dots \wedge du_n\,.\endgather$$
Here $\widetilde{p}(u_1,\dots ,u_n)$ denotes $p(1,u_1,\dots,u_n)$. We see that
the
form $\widehat{pr}^*\Phi^*\omega$ has the lo\-ga\-rithmic pole on the 
exceptional
divisor. The form ${\r_2}_{|\widehat{Y}\cup \Pt^n}$ is the second Leray 
residue; 
\cite{8}, \cite{9}.
We can decompose the form $\widehat{pr}^*\Phi^*\omega$ in a way
$$\widehat{pr}^*\Phi^*\omega=\frac{du_0}{u_0}\wedge
\frac{d\widetilde{\Phi^*s}}{\widetilde{\Phi^*s}}\wedge \r'_2\,,$$
where $\r'_2$ do not contain $u_0$ nor $du_0$. Then 
${\r'_2}_{|\widehat{Y}\cup
\Pt^n}= {\r_2}_{|\widehat{Y}\cup \Pt^n}$.
The function $\widetilde{\Phi^*(s)}$ describes 
$\widehat{Y}\cup \Pt^n$ in $\Pt^n$ for $u_0\neq 0$, so to show that
${r_2}_{|\widehat{Y}\cup \Pt^n} 
\not\equiv 0$ it suffices check that $d\widetilde{\Phi^*(s)}\wedge 
r'_2\not\equiv
0$ on $\widehat{Y}\cup \Pt^n$. By the decomposition:
$$u_0\,\widetilde{\Phi^*(s)}\,\widehat{pr}^*\Phi^*\omega=
du_0\wedge d\widetilde{\Phi^*s}\wedge \r'_2=
C\,du_0\wedge \widetilde{\Phi^*g}\left(\prod_{i=1}^{i=n}
u_i^{la_i-1}\right)\, du_1\wedge\dots\wedge du_n \,.$$
The polynomial $g$ has the lower weight than $s$, thus $\widetilde{\Phi^*g}$ 
does not vanish on
$\widehat{Y}\cup \Pt^n $. Moreover $\widehat{Y}\cup \Pt^n$ is not contained in
any of hyperplane $u_i=0$. Thus $d\widetilde{\Phi^*s}\wedge
\r'_2\not\equiv 0$ on 
$\widehat{Y}\cup \Pt^n $ and hence ${\r'_2}_{|\widehat{Y}\cup \Pt^n}\not\equiv
0$. This way we proved 

\proclaim{Theorem 7} Suppose that $g$ has the nonzero 
quasihomogeneous component of the
weight $\alpha=1-\kappa$. Then the second residue of $\omega=\frac
gsdz_0\wedge\dots\wedge dz_n$ does not vanish in $H^{n-1}(L/S^1)$.\endproclaim

\example{Example 8, \cite{13} }Consider the singularity of the type 
$P_8$: 
$$s(z_0,z_1,z_2)=z_0^3+z_1^3+z_2^3$$
and $\omega=\frac 1s dz_0\wedge dz_1\wedge dz_2$. Then $\widehat{L}/S^1=L/S^1
\subset \Pt^2$. The second residue (i.e. the obstruction to
lift) is:
$$res_2\omega=\left[\frac 1{2\pi i}\int_p Res\,\omega\right] =\frac
13\,u_1du_2-u_2du_1$$ 
in the notation used above. As one can check by hand the integral
$$\int_{L/S^1\cup \Rt\Pt^2}res_2\omega\neq 0\,.$$\endexample

\enddocument
\end